\pdfoutput=1
\RequirePackage{ifpdf}
\ifpdf % We are running pdfTeX in pdf mode
\documentclass[pdftex]{sigma}
\else
\documentclass{sigma}
\fi

\def\II{\hbox{{1}\kern-.25em\hbox{l}}}

\numberwithin{equation}{section}

\begin{document}

\allowdisplaybreaks

\newcommand{\arXivNumber}{1810.10806}

\renewcommand{\thefootnote}{}

\renewcommand{\PaperNumber}{121}

\FirstPageHeading

\ShortArticleName{Matrix Bailey Lemma and the Star-Triangle Relation}

\ArticleName{Matrix Bailey Lemma and the Star-Triangle Relation\footnote{This paper is a~contribution to the Special Issue on Elliptic Hypergeometric Functions and Their Applications. The full collection is available at \href{https://www.emis.de/journals/SIGMA/EHF2017.html}{https://www.emis.de/journals/SIGMA/EHF2017.html}}}

\Author{Kamil Yu.~MAGADOV~$^{\dag}$ and Vyacheslav P.~SPIRIDONOV~$^{\ddag\S}$}

\AuthorNameForHeading{K.Yu.~Magadov and V.P.~Spiridonov}

\Address{$^\dag$~Deceased; Moscow Institute of Physics and Technology,
Dolgoprudny, Moscow Region, Russia}

\Address{$^\ddag$~Laboratory of Theoretical Physics, JINR, Dubna, Moscow Region, 141980 Russia}
\EmailD{\href{mailto:spiridon@theor.jinr.ru}{spiridon@theor.jinr.ru}}

\Address{$^\S$~National Research University Higher School of Economics, Moscow, Russia}

\ArticleDates{Received August 10, 2018, in final form October 30, 2018; Published online November 10, 2018}

\Abstract{We compare previously found finite-dimensional matrix and integral operator
realizations of the Bailey lemma employing univariate
elliptic hypergeometric functions.
With the help of residue calculus we explicitly show how the integral
Bailey lemma can be reduced to its matrix version. As a consequence,
we demonstrate that the matrix Bailey lemma can be interpreted as a
star-triangle relation, or as a Coxeter relation for a permutation group.}


\Classification{33D60; 33E20}

\renewcommand{\thefootnote}{\arabic{footnote}}
\setcounter{footnote}{0}

\section{Introduction}

Quantum integrable systems can be realized as statistical mechanics models \cite{Baxter} and solved by the quantum version of the inverse scattering method \cite{TF}. One of the important examples is given by the hard hexagons model solved in \cite{Baxter2}. During its investigation Baxter met a particular type of $q$-hypergeometric series identities called the Rogers--Ramanujan identities~\cite{aar} without knowing how to prove them. It appeared that Bailey~\cite{Bailey1949} has found already a systematic way of proving such relations. The key ingredient in his method is the so-called Bailey lemma. Taking appropriate entries for this lemma, called Bailey pairs, Rogers--Ramanujan identities can be easily proved. An iterative mechanism for building new Bailey pairs from a given one has been found by Andrews \cite{Andrews} (see also~\cite{Paule}). Sequences of these pairs form the Bailey chain. A survey of the results derived from the Bailey lemma in the first 50 years was given by Warnaar in~\cite{Warnaar}.

Shortly after appearance of the summary \cite{Warnaar}, an extension of the Bailey chains technique to elliptic hypergeometric series was proposed by the second author in \cite{VS-IMRN-2002}. Some further analysis of the elliptic Bailey chains for such series was performed by Warnaar in~\cite{war:extensions}. A principal generalization of the whole Bailey formalism to the level of integrals and its application to derivation of infinite sequences of symmetry transformations for elliptic hypergeometric integrals was discovered by the second author in~\cite{spi:bailey}. It is curious to note that Bailey himself was interested in the generalization of his techniques from $q$-series to integrals, see the personal letters from Bailey to Dyson~\cite{Zudilin}. However, he was not able to find any use of such an idea, which was realized on the full scale only in \cite{spi:bailey}.

In \cite{DS} the integral Bailey pairs were used for building an integral operator (an $R$-matrix) satisfying the Yang--Baxter equation \cite{Baxter,TF}. The ingredients of the Bailey lemma serve as building blocks of this $R$-matrix. In this case the key operator identity generating an infinite sequence of Bailey pairs, known as the star-triangle relation, can be interpreted as an infinite-dimensional space realization of the cubic Coxeter relations for a permutation group. Such integral operator realizations of the generating relations of permutation groups were proposed earlier in \cite{D,Derkachov:2007gr} at the level of simpler special functions.

Constructions of the matrix and integral Bailey lemmas in \cite{VS-IMRN-2002} and \cite{spi:bailey} look similar in their structure. Nevertheless, the direct connection between them has not been described in detail yet. A residue calculus for a quite sophisticate type of (multivariate) integral Bailey lemmas was discussed by Rains in \cite{Rains}, but the exposition is not as explicit as one would wish to see. Our goal is to give an elementary description how an appropriate residue calculus applied to the integral Bailey lemma of \cite{spi:bailey} generates the matrix Bailey lemma of \cite{VS-IMRN-2002}. This gives an interpretation of the corresponding key matrix identity as a star-triangle relation which was not considered in~\cite{VS-IMRN-2002}. Consequently, we also interpret this relation as a realization of the Coxeter relations of a permutation group, similar to the integral case.

\section{Integral Bailey lemma and the star-triangle relation}

We start our consideration from the integral Bailey lemma derived in~\cite{spi:bailey}. For this we recall that two functions $\alpha(z,t)$ and $\beta(z,t)$ depending on two complex variables $z$ and $t$ are called the Bailey pair, if they are related by the following (univariate) integral transformation
\begin{gather}\label{EFT}
\beta(w,t)=M(t)_{wz}\alpha(z,t):=\kappa\int_\mathbb{T} \frac{\Gamma\big(tw^{\pm1}z^{\pm1};p,q\big)}{\Gamma\big(t^2,z^{\pm2};p,q\big)}\alpha(z,t)\frac{{\rm d}z}{z},\\
\kappa:=\frac{(p;p)_\infty(q;q)_\infty }{ 4\pi\textup{i} },\qquad (z;q)_\infty:=\prod_{j=0}^\infty\big(1-zq^j\big),\nonumber
\end{gather}
where $\mathbb{T}$ is the positive oriented unit circle, $|tw|, |t/w|<1$, and we assume that $\alpha(z,t)$ is analytical near $z\in\mathbb{T}$.

The kernel of the integral operator $M(t)_{wz}$ contains a combination of seven elliptic gamma functions
\begin{gather*}
\Gamma(z;p,q):=\prod_{j,k=0}^{\infty} \frac{1-z^{-1} p^{j+1}q^{k+1}}{1 - z p^j q^k}, \qquad |p|, |q|<1,
\end{gather*}
with the convention
\begin{gather*}
\Gamma(a,b;p,q):=\Gamma(a;p,q)\Gamma(b;p,q), \qquad \Gamma\big(tz^{\pm1};p,q\big):=\Gamma(tz;p,q)\Gamma\big(tz^{-1};p,q\big).
\end{gather*}
The integration contour $\mathbb{T}$ in \eqref{EFT} can be replaced by any contour $C$ encircling the same set of singularities of the kernel function which allows analytic continuation of the action of $M(t)_{wz}$ to a~wider range of values of parameters $t$ and $w$, provided $C$ does not cross singularities of $\alpha(z,t)$.

Let us list briefly the key properties of the elliptic gamma function. These are the symmetry in bases
\begin{gather*}
\Gamma(z;p,q)=\Gamma(z;q,p)
\end{gather*}
and the finite-difference equations
\begin{gather*}
\Gamma(qz;p,q)=\theta(z;p)\Gamma(z;p,q), \qquad\Gamma(pz;p,q)=\theta(z;q)\Gamma(z;p,q),%\label{eqs}
\end{gather*}
where $\theta(z;p)$ is the ``short'' Jacobi theta function
\begin{gather*}
 \theta(z;p):=(z;p)_\infty\big(pz^{-1};p\big)_\infty=\frac{1}{(p;p)_\infty} \sum_{n\in\mathbb{Z}}p^{n(n-1)/2}(-z)^n.%\label{Jacobi}
 \end{gather*}
Other properties we use include the inversion relation
\begin{gather}
\Gamma(z;p,q)=\frac{1}{\Gamma\big(\frac{pq}{z};p,q\big)}, \label{inv}
\end{gather}
the quadratic transformation (note its appearance in \eqref{EFT})
\begin{gather*}
\Gamma\big(z^2;p,q\big)=\Gamma\big(\pm z,\pm q^{1/2}z,\pm p^{1/2}z,\pm (pq)^{1/2}z;p,q\big),
\end{gather*}
and the limiting relation
\begin{gather*}
\lim_{z\to1}(1-z)\Gamma(z;p,q)=\frac{1}{(p;p)_\infty (q;q)_\infty},
\end{gather*}
which is needed in the residue calculus.

The integral transformation \eqref{EFT} is also called the elliptic Fourier transformation. One of the justification for such a name comes from the fact established in~\cite{spi-war:inversions} that for actions on the $A_1$-symmetric functions, $f\big(z^{-1}\big)=f(z)$, and under some constraints on the parameters, the inversion of this transform is obtained by the reflection $t\to t^{-1}$, i.e.,
\begin{gather*}
M(t)_{wz}^{-1}=M\big(t^{-1}\big)_{wz},\qquad \text{or}\qquad M(t)_{wz}M\big(t^{-1}\big)_{zx}f(x)=f(w),
\end{gather*}
where we assume an analytic continuation of the action of the $M$-operators by an appropriate deformation of the contours of integration and test functions $f(x)$ with relevant analytical properties. For a discussion of this procedure and partial description of the null space of the integral operator~$M(t)_{wz}$, see~\cite{DStmf}.

Let us define now the following combination of four elliptic gamma functions
\begin{gather*}
D(s;y,w):=\Gamma\big(\sqrt{pq}s^{-1}y^{\pm{1}}w^{\pm{1}};p,q\big).%\label{D}
\end{gather*}
From the inversion relation \eqref{inv} we find that
\begin{gather*}
D\big(s^{-1};y,w\big)=\frac{1}{D(s;y,w)}.
\end{gather*}

The integral Bailey lemma established in \cite{spi:bailey} states that from a given Bailey pair $\alpha(w,t)$ and $\beta(w,t)$ one can obtain a new Bailey pair with the replacement of the parameter $t$ by another arbitrary variable. Namely, the functions
\begin{gather}\label{B1}
\alpha'(w,st):=D(s;y,w)\alpha(w,t),\\
\beta'(w,st):=D\big(t^{-1};y,w\big)M(s)_{wx}D(st;y,x)\beta(x,t),\label{B2}
\end{gather}
with the assumption $|x|=1$ and the constraints $|t|, |s|, \big|\sqrt{pq}s^{-1}t^{-1}y^{\pm{1}}\big| <1$, form a new Bailey pair, i.e.,
\begin{gather*}
\beta'(w,st)=M(st)_{wz}\alpha'(z,st).
\end{gather*}
The transformed functions $\alpha'$ and $\beta'$ depend now on two new complex variables $s$ and $y$. Evidently, this procedure has an iterative character, i.e., from a given Bailey pair it generates infinitely many such pairs containing as many free variables as many times we apply the maps \eqref{B1} and \eqref{B2}. Note that for taken restrictions on parameters one automatically has $\big|\sqrt{pq}s^{-1}y^{\pm{1}}\big|<1$ and the unit circle separates sequences of poles of $D(st;y,x)$ and $D(s;y,x)$ converging to $x=0$ from their reciprocals.

It is easy to see that this statement leads to the following operator identity
\begin{gather} \label{Bailey}
M(s)_{wx}D(st;y,x)M(t)_{xz}=D(t;y,w)M(st)_{wz}D(s;y,z),
\end{gather}
which is called the star-triangle relation. The proof of \eqref{Bailey} is based on the explicit elliptic beta integral evaluation formula derived in \cite{VS-2001}. Namely, for six complex parameters $t_1,t_2,\dots,t_6$ subject to the constraints $\prod\limits_{j=1}^6 t_j=pq$ and $|t_j|<1$ one has
\begin{gather}\label{EBT2}
\kappa\int_{\mathbb{T}}\frac{\prod\limits_{j=1}^6\Gamma\big(t_jz^{\pm{1}};p,q\big)}{\Gamma\big(z^{\pm{2}};p,q\big)} \frac{{\rm d}z}{z}=\prod_{1\leq j<k\leq 6}\Gamma(t_jt_k;p,q).
\end{gather}
In a special limit $p\rightarrow 0$ formula \eqref{EBT2} is reduced to the Rahman $q$-beta integral~\cite{rah:integral}.

In \cite{DS} the star-triangle relation \eqref{Bailey} was used for a rigorous construction of an integral operator solving the Yang--Baxter equation, which is currently the most complicated known rank 1 solution of this important equation of mathematical physics. Earlier the functional form of relation~\eqref{Bailey} was applied for building new integrable two-dimensional lattice spin systems, see survey \cite{BKS} and references therein.

Following \cite{DS}, let us recall how equality \eqref{Bailey} can be interpreted as the cubic Coxeter relation for the permutation group $\mathbb{S}_3$. Consider three complex variables $\underline{u}=(u_1,u_2,u_3)$ and define the elementary generators of the permutation group
\begin{gather*}
s_1\underline{u}=(u_2,u_1,u_3), \qquad s_2\underline{u}=(u_1,u_3,u_2).
\end{gather*}
Let us denote
\begin{gather*}
s={\rm e}^{2\pi \textup{i}a}, \qquad t={\rm e}^{2\pi \textup{i}b}, \qquad a=u_2-u_3, \qquad b=u_1-u_2.
\end{gather*}
Now one can introduce the integral operators $S_1$ and $S_2$, acting as
\begin{gather*}%\label{Def_Kokster}
 [S_1(\underline{u})f](w):=M(t)_{wx}f(x),\qquad [S_2(\underline{u})f](w):=D(s;y,w)f(w).
\end{gather*}
Here one can write symbolically $S_2$ as an integral operator with a singular (delta-function) kernel, $S_2(\underline{u})_{wx}=D(s;y,x)\II_{wx}$, where $\II:=\II_{wx}$ is the unit operator acting in the infinite-dimensional
space of complex functions, $\II_{wx}f(x)=f(w)$.

The operators $S_1(\underline{u})$ and $S_2(\underline{u})$ depend not on all variables $u_j$, but on their particular combinations, namely:
\begin{gather*}%\label{Def2Cox}
S_1(\underline{u})=S_1(u_1-u_2), \qquad S_2(\underline{u})=S_2(u_2-u_3).
\end{gather*}
One takes the following multiplication rule for these operators~\cite{Derkachov:2007gr,DS}:
\begin{gather*}%\label{Multiply}
S_jS_k:=S_j(s_k\underline{u})S_k(\underline{u}), \qquad j, k=1,2.
\end{gather*}
Then it is easy to see that
\begin{gather*}
S_2^2:=S_2(-a)S_2(a)=\II, \qquad a=u_2-u_3.
\end{gather*}

As mentioned above, for functions satisfying the restriction $f\big(z^{-1}\big)=f(z)$ and under some constraints on the values of $a$ and the contours of integration \cite{DStmf, spi-war:inversions}, one has also the inversion relation
\begin{gather*}
S_1^2:=S_1(-b)S_1(b)=\II, \qquad b=u_1-u_2.
\end{gather*}
Using the accepted multiplication rule, the relation \eqref{Bailey} can be written as follows \cite{DS}
\begin{gather*}
S_1S_2S_1:=S_1(s_2s_1\underline{u})S_2(s_1\underline{u})S_1(\underline{u})=S_1(a)S_2(a+b)S_1(b)\\
\hphantom{S_1S_2S_1}{} =S_2(b)S_1(a+b)S_2(a)=S_2(s_1s_2\underline{u})S_1(s_2\underline{u})S_2(\underline{u})=:S_2S_1S_2.%\label{Cox=STR}
\end{gather*}
So, the integral Bailey lemma operators $D$ and $M$ determine the generators of elementary permutations of parameters acting as integral operators in the infinite-dimensional space of complex functions, and the basic relation~\eqref{Bailey} is equivalent to the cubic Coxeter relation $S_1S_2S_1=S_2S_1S_2$.

\section{Matrix reduction of the integral Bailey pairs}

We consider possible reductions of the integral operator $M(t)_{xz}$ to simpler matrix forms. For that we discuss first analytic continuation of the action of $M(t)_{xz}$. The kernel of this integral operator has poles in the integration variable at the following points
\begin{gather*}
z=\big\{tx^{\pm1}q^ap^b, t^{-1}x^{\pm1}q^{-a}p^{-b}\big\}_{a,b\in \mathbb{Z}_{\geq 0}}.
\end{gather*}
Using the inversion formula for the $\Gamma$-function, one finds
\begin{gather*}
\frac{1}{\Gamma\big(z^{\pm2};p,q\big)}=\theta\big(z^2;q\big)\theta\big(z^{-2};p\big).
\end{gather*}
Therefore zeroes of the kernel are located at the points:
\begin{gather*}
z=\big\{\pm q^{a/2}, \pm p^{a/2}\big\}_{a\in \mathbb{Z}}\cup \big\{ t^{-1}x^{\pm 1}q^{a+1}p^{b+1}, tx^{\pm 1}q^{-a-1}p^{-b-1} \big\}_{a,b\in \mathbb{Z}_{\geq 0}}.
\end{gather*}
Evidently, at the point $z=0$ there is an essential singularity. In the definition~\eqref{EFT} we assumed that $\big|tw^{\pm{1}}\big|<1$, in which case the contour $\mathbb{T}$ separates all geometric sequences of poles accumulating to zero from their reciprocals that go to infinity. We can deform $\mathbb{T}$ to different contours as long as we do not cross poles of the kernel and integrated function and this provides the desired analytic continuation
of the action of $M(t)_{xz}$ to wider ranges of the variables~$t$ and~$x$.

Now we assume that the function $\alpha(z)$ in \eqref{EFT} satisfies the condition $\alpha(z)=\alpha\big(z^{-1}\big)$. The reason for such a restriction comes from the fact that the image of operator $M(t)_{xz}$ obeys this symmetry. Moreover, the inversion relation $M^{-1}(t)_{xz}=M\big(t^{-1}\big)_{xz}$ holds only for such functions.

Let us suppose that the poles of $\alpha(z)$ are simple and that they are located at the points $\big\{ z=z_{0}q^{m} \big\}_{m=0,1,\dots,N}$ with $|z_0|<1$ and their $z\to 1/z$ reciprocals. Now we deform the contour of integration $\mathbb{T}$ to a contour $C$ that crosses these poles of $\alpha(z)$ without touching any other pole. From the Cauchy theorem we obtain the equality
 \begin{gather}
\beta(x)= M(t)_{xz}\alpha(z)=\kappa\int_{C}\frac{\Gamma\big(tx^{\pm{1}}z^{\pm{1}};p,q\big)}{\Gamma\big(t^{2},z^{\pm{2}};p,q\big)}\alpha(z)\frac{{\rm d}z}{z}\nonumber\\
\hphantom{\beta(x)=}{}+4\pi \textup{i}\kappa\sum_{m=0}^{N}\lim_{z\to z_{0}q^{m}}\big(1-z_{0}q^{m}/z\big)\frac{\Gamma\big(tx^{\pm{1}}z^{\pm{1}};p,q\big)}{\Gamma\big(t^{2},z^{\pm{2}};p,q\big)}\alpha(z)\nonumber\\
\hphantom{\beta(x)}{}=\kappa\int_{C}\frac{\Gamma\big(tx^{\pm{1}}z^{\pm{1}};p,q\big)}{\Gamma\big(t^{2},z^{\pm{2}};p,q\big)}\alpha(z)\frac{{\rm d}z}{z} +4\pi\textup{i}\kappa\sum_{m=0}^{N}
\frac{\Gamma\big(tx^{\pm1}(z_0q^m)^{\pm1};p,q\big)}{\Gamma\big(t^2,(z_0q^m)^{\pm2};p,q\big)}\tilde\alpha_m,\label{Res1}
\end{gather}
where $\tilde\alpha_m=\lim\limits_{z\to z_{0}q^{m}}\big(1-z_{0}q^{m}/z\big)\alpha(z)$ are the residues of poles of $\alpha(z)/z$.

\begin{figure}[t]\centering
\includegraphics[width=10cm]{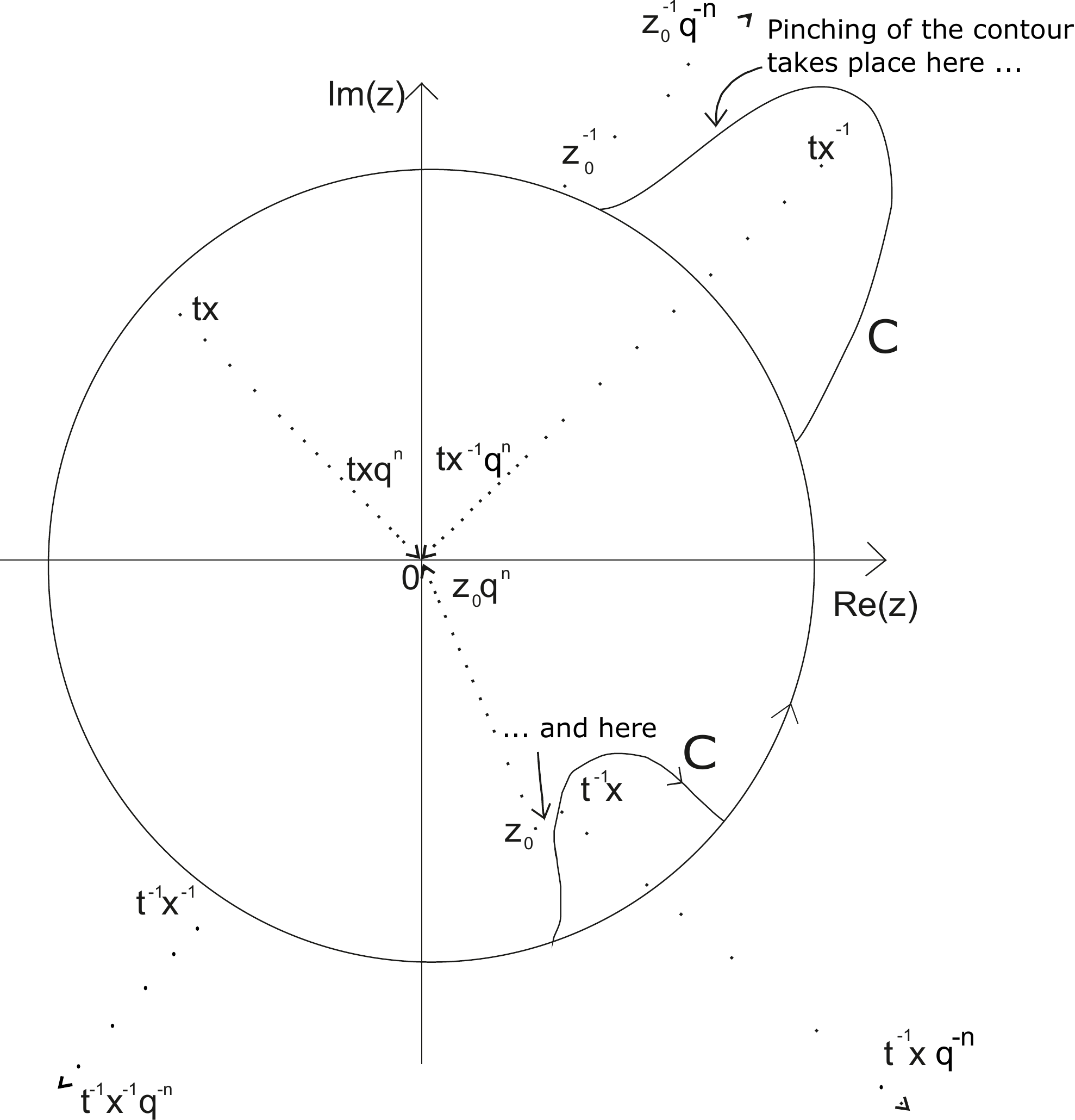}
\caption{Pinching of the integration contour $C$ by sequences of poles $\big\{z=tx^{-1}q^n, z_0^{-1}q^{-m}\big\}$ and $\big\{z=t^{-1}xq^{-n}, z_0q^m\big\}$ in the limit $\epsilon\to 1$.}\label{figure}
\end{figure}

In order to match with the notation of \cite{VS-IMRN-2002} we set
\begin{gather*}
z_0=a^{1/2}, \qquad t=\frac{k^{1/2}}{a^{1/2}}, \qquad x=\epsilon^{-1} k^{1/2} q^N,
\end{gather*}
where $N$ is the same integer number as in \eqref{Res1}. The sum of residues takes the following form
\begin{gather*}
4\pi\textup{i}\kappa\sum_{m=0}^{N}\frac{\Gamma\big(\epsilon q^{m-N},k\epsilon^{-1}q^{m+N},a^{-1}\epsilon q^{-m-N},ka^{-1}\epsilon^{-1} q^{N-m};p,q\big)}{\Gamma\big(k/a,\big(aq^{2m}\big)^{\pm1};p,q\big)}\tilde\alpha_m.
\end{gather*}
For $\epsilon$ sufficiently close to 1 the condition $|tx^{-1}|<1$ will be broken and in the definition of our $M$-operator we have to assume that the contour $C$ does not cross corresponding poles lying outside~$\mathbb{T}$. However, in the limit $\epsilon\to 1$ we have the following geometrical picture: before we force the contour $C$ to cross poles of $\alpha(z)$ it gets pinched by $2(N+1)$ poles of $\alpha(z)$ from one side (the doubling comes from the symmetry $z\to1/z$) and $2(N+1)$ poles of the kernel of $M$-operator from the other side, as described on the Fig.~\ref{figure} for the choice $0<q<1$.

The multiplier $\Gamma\big(tx^{-1}z_0q^m\big)=\Gamma\big(\epsilon q^{m-N}\big)$ in the residue sum is diverging in the limit $\epsilon\to1$. To remove this divergency we multiply the above expression by $1-\epsilon$ and take the needed limit. Denote
\begin{gather*}
\tilde\beta_N(a,k):= \lim_{\epsilon\to1}(1-\epsilon)\beta(x)
\end{gather*}
and use the relation
\begin{gather*}
\lim_{\epsilon\to1}(1-\epsilon)\Gamma\big(\epsilon q^{m-N};p,q\big) =\frac{1}{\theta\big(q^{m-N},q^{m-N+1},\ldots, q^{-1};p\big)(q;q)_\infty(p;p)_\infty},
\end{gather*}
with the convention $\theta(a,b;p):=\theta(a;p)\theta(b;p)$. Since the integral part in~\eqref{Res1} is finite for $\epsilon\to1$, multiplication by $1-\epsilon$ and subsequent limit $\epsilon\to1$ kills it completely. As a result we obtain
\begin{gather} \label{Res2}
\tilde{\beta}_N(a,k)=\sum_{m=0}^N \frac{\Gamma\big(kq^{N+m};p,q\big)\Gamma\big(\frac{k}{a}q^{N-m};p,q\big)\Gamma\big(a^{-1}q^{-N-m};p,q\big)}{\theta\big(q^{m-N},q^{m-N+1},\ldots, q^{-1};p\big)\Gamma(k/a;p,q)\Gamma\big(\big(aq^{2m}\big)^{\pm1};p,q\big)}\tilde{\alpha}_{m}(a,k).
\end{gather}

Now we apply the inversion relation \eqref{inv}, the equation
\begin{gather*}
\Gamma\big(zq^{n};p,q\big)=\theta(z;p)_{n}\Gamma(z;p,q),
\end{gather*}
where
\begin{gather*}
\theta(z;p)_n =\theta(z)_n=
\begin{cases}
\displaystyle \prod_{j=0}^{n-1}\theta\big(zq^j;p\big), & n>0, \\ \displaystyle
\prod_{j=1}^{-n}\frac{1}{\theta\big(zq^{-j};p\big)}, & n<0
\end{cases}
\end{gather*}
is the elliptic Pochhammer symbol, and use the quasiperiodicity properties of $\theta(z;q)$-function. After corresponding simplifications, relation~(\ref{Res2}) reduces to the following form
\begin{gather*}
\tilde{\beta}_{N}(a,k)=\frac{\Gamma(k;p,q)} {\Gamma(a;p,q)}\sum_{m=0}^{N}M_{Nm}(a,k)q^{N(N+1)-m(m-1)}\tilde{\alpha}_{m}(a,k),
\end{gather*}
where
\begin{gather}\label{matrix}
M_{Nm}(a,k)=\frac{\theta(k)_{N+m}\theta\big(\frac{k}{a}\big)_{N-m}}{\theta(qa)_{N+m}\theta(q)_{N-m}}\frac{\theta\big(aq^{2m};p\big)}{\theta(a;p)}a^{N-m}.
\end{gather}
Let us renormalize the sequences $\tilde\alpha_n$ and $\tilde\beta_n$,
\begin{gather*}
\beta_{N}(a,k):=\tilde{\beta}_{N}(a,k)q^{-N(N+1)}, \qquad \alpha_{m}(a,k):= \tilde{\alpha}_{m}(a,k) q^{-m(m+1)}\frac{\Gamma(k;p,q)}{\Gamma(a;p,q)}.
\end{gather*}
Finally, we obtain the matrix relation
\begin{gather}\label{DBP}
\beta_{N}(a,k)=\sum_{m=0}^{N}M_{Nm}(a,k)\alpha_{m}(a,k).
\end{gather}
This is precisely the definition of the discrete Bailey pairs $(\alpha_n,\beta_n)$ introduced in~\cite{VS-IMRN-2002}, where the matrix $M_{Nm}(a,k)$ has appeared for the first time. For $p\to 0$ it reduces to the Bressoud matrix~\cite{bre}. Note that due to the presence of the factor $1/\theta(q)_{N-m}$ matrix elements~\eqref{matrix} vanish for $m>N$, which reduces the matrix to a triangular form.

Thus we have shown how in the appropriate limit the elliptic Fourier transformation operator is reduced to a matrix action (a discrete finite-difference operator). A different type of reduction of the $M$-operator has been considered in~\cite{DS}. Namely, in the limit $t^2\to q^{-N}$, $N\in \mathbb{Z}_{\geq 0}$, sequences of poles $z=txq^n$ and $z=t^{-1}xq^{-m}$ and their reciprocals start to pinch the contour of integration and there emerges also the vanishing factor $1/\Gamma(t^2;p,q)$. Forcing the latter contour to cross half of these poles and taking the limit one finds
\begin{gather*}
[M(t)f](x) =\frac{\Gamma\big(x^{-2};p,q\big)}{\Gamma\big(t^{-2}x^{-2};p,q\big)}\sum_{k=0}^{N}\frac{\theta\big((tx)^2q^{2k};p\big)} {\theta\big((tx)^2;p\big)}\frac{\theta\big(t^2,(tx)^2;p\big)_k} {\theta\big(q,qx^2;p\big)_k}
\frac{f\big(tq^kx\big)}{t^{4k}x^{2k}q^{k^2}},
\end{gather*}
where $t=\pm q^{-N/2}$. Note that this is an analytical finite-difference operator acting on functions of complex variable $x$, whereas in our present case we have a matrix acting on sequences of numbers.

In~\cite{DS} an even more complicated reduction associated with the limit $t^2\to q^{-N}p^{-M}$, $N, M\in \mathbb{Z}_{\geq 0}$, was considered. It results in the finite-difference operator involving scalings by two independent variables $q$ and $p$ which we do not present here. We only note that a similar complicated two-base reduction can be applied to the matrix reductions considered in this note. This should lead to a ``doubling'' of our matrices, but we skip such an opportunity for now.

\section{Matrix reduction of the integral Bailey lemma}

Denoting as $\alpha$ and $\beta$ columns formed by $\alpha_n$ and $\beta_n$, $n=0,\ldots, N$, we can rewrite relation \eqref{DBP} in the matrix form $\beta(a,k)=M(a,k)\alpha(a,k)$. Define the diagonal matrix
\begin{gather*}
D_{nm}(a;b,c)=D_m(a;b,c)\delta_{nm}, \qquad D_m(a;b,c)=\frac{\theta(b,c)_m}{\theta(aq/b,aq/c)_m} \left(\frac{aq}{bc}\right)^m.
\end{gather*}
The discrete Bailey lemma of \cite{VS-IMRN-2002} states that from a given Bailey pair $\alpha(a,t)$ and $\beta(a,t)$ one can form infinitely many such pairs using the following transformation rules
\begin{gather*}%\label{rel1}
 \alpha'(a,k)=D(a;b,c)\alpha(a,t),\\
\beta'(a,k)=D(k;qt/b,qt/c)M(t,k)D(t;b,c)\beta(a,t),
\end{gather*}
where $k$, $b$, $c$ are arbitrary new parameters satisfying the constraint $kbc=qat$. Note that the parameter $t$ in the Bailey pairs is replaced by a new variable $k$. This means validity of the relations $\beta'(a,k)=M(a,k)\alpha'(a,k)$, which leads to the following
matrix identity
\begin{gather}\label{MatrixBailey}
M(a,k)D(a;b,c)M(t,a)=D(k;qt/b,qt/c)M(t,k)D(t;b,c).
\end{gather}
After substitution of the explicit expressions for matrices, one can see that it holds true due to the Frenkel--Turaev summation formula~\cite{ft}.

Let us reduce now relation (\ref{Bailey}) and show how it generates identity \eqref{MatrixBailey}. We consider the left-hand side of equality \eqref{Bailey} and deform the contour of integration for $z$-variable $\mathbb{T}$ to~$C_z$ in the same way as we did before by crossing the poles of $\alpha(z)$ at $z=z_0q^m$, $m=0,\ldots, N,$ and their reciprocals. For simplicity of the analysis of the structure of pole sequences we will assume that $0<q<1$. In the final results it will be easy to make analytic continuation to arbitrary complex~$q$. After the residue calculus we obtain
\begin{gather}
 M(s)_{wx}D(st;y,x)M(t)_{xz}\alpha(z)=4\pi \textup{i}\kappa^{2}\int_{\mathbb{T} }\frac{{\rm d}x}{x}\frac{\Gamma\big(sw^{\pm{1}}x^{\pm{1}};p,q\big)}
{\Gamma\big(s^2,x^{\pm{2}};p,q\big)}\Gamma\big(\sqrt{pq}s^{-1}t^{-1}y^{\pm{1}}x^{\pm{1}};p,q\big)\nonumber \\
\qquad{} \times \sum_{m=0}^N\tilde{\alpha}_{m} \frac{\Gamma\big(tx^{\pm{1}}\big(z_{0}q^{m}\big)^{\pm{1}};p,q\big)}{\Gamma\big(t^{2},\big(z_{0}q^{m}\big)^{\pm{2}};p,q\big)}
+{\kappa}^{2}\int_{\mathbb{T} }\frac{{\rm d}x}{x}\frac{\Gamma\big(sw^{\pm{1}}x^{\pm{1}};p,q\big)}{\Gamma\big(s^2,x^{\pm{2}};p,q\big)}\nonumber\\
\qquad{} \times \Gamma\big(\sqrt{pq}s^{-1}t^{-1}y^{\pm{1}}x^{\pm{1}};p,q\big)\int_{C_z}\frac{{\rm d}z}{z}\frac{\Gamma\big(tx^{\pm{1}}z^{\pm{1}};p,q\big)}{\Gamma\big(t^{2},z^{\pm{2}};p,q\big)}\alpha(z),
\label{redMDM1} \end{gather}
where we denoted as before $\tilde{\alpha}_{m}=\lim\limits_{z\to z_{0}q^{m}}\big(1-z_{0}q^{m}/z\big)\alpha(z)$.

Consider the first term of this expression. Its integrand has poles at the points
\begin{gather*}
x=sw^{\pm1}p^jq^k,\qquad \sqrt{pq}s^{-1}t^{-1}y^{\pm1}p^jq^k,\qquad t\big(z_0q^m\big)^{\pm1}p^jq^k,\qquad j,k\in\mathbb{Z}_{\geq0},
\end{gather*}
and their $x\to 1/x$ reciprocals. Now we deform the $x$-integration contour from $\mathbb{T} $ to $C_x$ which crosses the poles $x=tz_0q^n$, $n=m,m+1,\ldots, N$, and their $x\to 1/x$ reciprocals without touching any other singularity. Again applying the residue calculus we come to the expression
\begin{gather}\nonumber
 M(s)_{wx}D(st;y,x)M(t)_{xz}\alpha(z)\\ \nonumber
{} =(4\pi \textup{i}\kappa)^2\sum_{m=0}^N\sum_{n=m}^N\tilde{\alpha}_m\lim_{x\to tz_{0}q^{n}}\left(1-\frac{tz_{0}q^n}{x}\right)\Gamma\left(\frac{tz_{0}q^m}{x};p,q\right)
\frac{\Gamma\big(sw^{\pm{1}}x^{\pm{1}};p,q\big)}{\Gamma\big(s^2,x^{\pm{2}};p,q\big)}\\
 {} \times \Gamma\big(\sqrt{pq}s^{-1}t^{-1}y^{\pm{1}}x^{\pm{1}};p,q\big)\frac{\Gamma\big(tx\big(z_{0}q^{m}\big)^{\pm{1}};p,q\big)\Gamma\big(tx^{-1}\big(z_{0}q^{m}\big)^{-1};p,q\big)}
{\Gamma\big(t^{2},\big(z_{0}q^{m}\big)^{\pm{2}};p,q\big)}\label{redMDM} \\
{}+ 4\pi \textup{i}\kappa^2\int_{C_x}\frac{{\rm d}x}{x} \frac{\Gamma\big(sw^{\pm{1}}x^{\pm{1}};p,q\big)}{\Gamma\big(s^2,x^{\pm{2}};p,q\big)}\Gamma\big(\sqrt{pq}s^{-1}t^{-1}y^{\pm{1}}x^{\pm{1}};p,q\big)\nonumber
\sum_{m=0}^{N}\tilde{\alpha}_m\frac{\Gamma\big(tx^{\pm{1}}\big(z_{0}q^{m}\big)^{\pm{1}};p,q\big)}{\Gamma\big(t^{2},\big(z_{0}q^{m}\big)^{\pm{2}};p,q\big)}\\ \nonumber
+{\kappa}^{2}\int_{\mathbb{T} }\frac{{\rm d}x}{x}\frac{\Gamma\big(sw^{\pm{1}}x^{\pm{1}};p,q\big)}{\Gamma\big(s^2,x^{\pm{2}};p,q\big)}\Gamma\big(\sqrt{pq}s^{-1}t^{-1}y^{\pm{1}}x^{\pm{1}};p,q\big)
\int_{C_z}\frac{{\rm d}z}{z}\frac{\Gamma\big(tx^{\pm{1}}z^{\pm{1}};p,q\big)}{\Gamma\big(t^{2},z^{\pm{2}};p,q\big)}\alpha(z).
\end{gather}

Using the relations
\begin{gather*}
\lim_{x\to tz_{0}q^{n}}\left(1-\frac{tz_{0}q^n}{x}\right)\Gamma\left(\frac{tz_{0}q^m}{x};p,q\right)=\frac{1}{(p;p)_{\infty}(q;q)_{\infty} \theta\big(q^{m-n}\big)_{n-m}},\\
\theta\big(q^{m-n}\big)_{n-m}=(-1)^{n-m}q^{-\frac{(n-m)(n-m+1)}{2}}\theta(q)_{n-m},
\end{gather*}
one can see that there is the multiplier $1/\theta(q)_{n-m}$ in the sum which vanishes for $n<m$. Therefore the summation over $n$ actually starts from zero and we can interchange summations over $m$ and $n$.

In the reduction of Bailey pairs considered in the previous section we introduced an $\epsilon$-parametrization of the external coordinate of the $M$-operator with the subsequent limit $\epsilon\to 1$, which converted the integral operator to a matrix. In the present case it works as follows. Namely, we fix
\begin{gather*}
w= \epsilon^{-1}stz_0q^N, \qquad \text{i.e.}, \qquad \Gamma\big(sw^{-1}tz_0q^n;p,q\big)=\Gamma\big(\epsilon q^{-N+n};p,q\big),
\end{gather*}
multiply the whole expression by $1-\epsilon$, and take the limit $\epsilon\to 1$. It can be seen that during the deformation of $x$-variable integration contour $\mathbb{T}$ to $C_x$ in passing from~\eqref{redMDM1} to~\eqref{redMDM} the limit $\epsilon \to 1$ leads to pinching of this contour precisely as it took place in the previous section (see Fig.~\ref{figure}), i.e., computation of the residues at $x=tz_0q^n$ is inevitable. In the same way as before it shows that the second term in~\eqref{redMDM} with the integral over $x$ remains finite for $\epsilon\to 1$, because~$C_x$ crossed already dangerous poles without touching other singularities. Similar situation holds with the third term in~\eqref{redMDM} containing integrals over~$x$ and~$z$. Namely, since $\big|sw^{-1}\big|$ becomes bigger than~1, we have to deform $x$-integration contour appropriately not to cross over any pole. This is possible due to the choice of the $z$-variable integration contour~$C_z$, which escapes the dangerous $z$-values $z=\big(z_0q^m\big)^{\pm1}$, and a special choice of the function~$\alpha(z)$. The constraint $0<q<1$ helps to make such a statement transparent. Therefore after multiplication by $1-\epsilon$ both integral terms in~\eqref{redMDM} vanish in the limit $\epsilon\to 1$.

As a result, we find the following replacements of the integral operator factors
\begin{gather*}
M(t)_{xz} \to \frac{\Gamma\big(t^2z_0^2;p,q\big)}{\Gamma\big(z_0^2;p,q\big)}\sum_{m=0}^N M_{nm}\big(z_0^2,t^2z_0^2\big)q^{n(n+1)-m(m+1)},\\
M(s)_{wx}\to \frac{\Gamma\big(s^2x_0^2;p,q\big)}{\Gamma\big(x_0^2;p,q\big)}\sum_{n=0}^N M_{Nn}\big(x_0^2,s^2x_0^2\big)q^{N(N+1)-n(n+1)},
\end{gather*}
where $x_0=tz_0$. For reducing the $D$-function multiplier, we have
\begin{gather*}
D(st;y,x)\to \Gamma\big(\sqrt{pq}s^{-1}t^{-1}y^{\pm 1}\big(x_0q^n\big)^{\pm1};p,q\big).
\end{gather*}
After denoting $a:=x_0^2$, $k:=s^2x_0^2$, $\tilde t:=z_0^2$ we obtain
\begin{gather*}
D(st;y,x)\to \frac{\Gamma(b,pc;p,q)}{\Gamma(qa/c, pqa/b;p,q)} D_{n}(a;b,c),\\
D_{n}(a;b,c):=\frac{\theta(b)_{n}\theta(c)_{n}}{\theta(qa/c)_{n}\theta(qa/b)_{n}}\left(\frac{qa}{bc}\right)^{n},
\end{gather*}
where we have introduced two parameters $b$ and $c$,
\begin{gather}\label{bc}
b:=\sqrt{\frac{pq\tilde t a}{k}}y, \qquad c:=\sqrt{\frac{q\tilde ta}{pk}} y^{-1},
\end{gather}
satisfying the important product rule
\begin{gather*}
kbc=qa\tilde t.
\end{gather*}

The total expression for the left-hand side takes the following form
\begin{gather}\nonumber
 M(s)_{wx}D(st;y,x)M(t)_{xz}\alpha(z) \\
 \qquad{} \to \frac{\Gamma(k,b,pc;p,q)q^{N(N+1)}}{\Gamma\big(\tilde t,qa/c,pqa/b;p,q\big)}\sum_{n,m=0}^N M_{Nn}(a,k)D_{n}(a;b,c)M_{nm}\big(\tilde t,a\big)q^{-m(m+1)}\tilde\alpha_m.\label{LHS}
 \end{gather}

Now we consider the right-hand side of equality (\ref{Bailey}). Similarly to the previous case, we deform the integration contour of $M(st)_{wz}$ from $\mathbb{T}$ to such a contour~$C$ that crosses the poles $z=z_0q^m$, $m=0,\ldots, N,$ of $\alpha(z)$ and does not touch other singularities. As a result, we obtain
\begin{gather}
D(t;y,w)M(st)_{wz}D(s;y,z)\alpha(z)=\kappa\int_{\mathbb{T}} \Gamma\big(\sqrt{pq}t^{-1}y^{\pm{1}}w^{\pm{1}};p,q\big)\nonumber \\
\qquad{}\times \frac{\Gamma\big(stw^{\pm{1}}z^{\pm{1}};p,q\big)}{\Gamma\big(s^2t^2,z^{\pm{2}};p,q\big)}\Gamma\big(\sqrt{pq}s^{-1}y^{\pm{1}}z^{\pm{1}};p,q\big)\alpha(z)\frac{{\rm d}z}{z}=\kappa\int_{C}
\Gamma\big(\sqrt{pq}t^{-1}y^{\pm{1}}w^{\pm{1}};p,q\big)\nonumber \\
\qquad{}\times \frac{\Gamma\big(stw^{\pm{1}}z^{\pm{1}};p,q\big)}{\Gamma\big(s^2t^2,z^{\pm{2}};p,q\big)}\Gamma\big(\sqrt{pq}s^{-1}y^{\pm{1}}z^{\pm{1}};p,q\big)\alpha(z)\frac{{\rm d}z}{z}
+\Gamma\big(\sqrt{pq}t^{-1}y^{\pm{1}}w^{\pm{1}};p,q\big)\nonumber \\
\qquad{}\times 4\pi\textup{i}\kappa\sum_{m=0}^{N}\tilde{\alpha}_m \frac{\Gamma\big(stw^{\pm{1}}\big(z_{0}q^{m}\big)^{\pm{1}};p,q\big)}{\Gamma\big(s^2t^2,\big(z_{0}q^{m}\big)^{\pm{2}};p,q\big)}
\Gamma\big(\sqrt{pq}s^{-1}y^{\pm{1}}\big(z_{0}q^{m}\big)^{\pm{1}};p,q\big).\label{RHS}
\end{gather}
We set as before $w= \epsilon^{-1}stz_0q^N$, multiply the whole expression by $1-\epsilon$ and take the limit $\epsilon\to 1$. Again the integral part in~\eqref{RHS} is finite and disappears in this limit. After detailed computations, we find the following replacements of the operator factors
\begin{gather*}
D(t;y,w)\to \frac{\Gamma\big(q\tilde t/c,pq\tilde t/b;p,q\big)}{\Gamma\big(kb/\tilde t, pkc/\tilde t;p,q\big)} D_{N}\big(k;q\tilde t/c,q\tilde t/b\big),\\
M(st)_{wz} \to \frac{\Gamma(k;p,q)}{\Gamma\big(\tilde t;p,q\big)}\sum_{m=0}^N M_{Nm}\big(\tilde t,k\big)q^{N(N+1)-m(m+1)},\\
D(s;y,z) \to \frac{\Gamma(b,pc;p,q)}{\Gamma\big(q\tilde t/c, pq\tilde t/b;p,q\big)} D_{m}(\tilde t;b,c).
\end{gather*}
The whole right-hand side expression takes the form
\begin{gather}
\lim_{\epsilon \to 1}(1-\epsilon)\, { \rm r.h.s.\, expression} = \sum_{m=0}^{N}\frac{\Gamma(k,b,pc;p,q)}{\Gamma\big(\tilde{t},\frac{qa}{c},\frac{pqa}{b};p,q\big)}q^{N(N+1)-m(m+1)}\nonumber\\
\hphantom{\lim_{\epsilon \to 1}(1-\epsilon)\, { \rm r.h.s.\, expression} =}{} \times D_{N }\big(k;q\tilde t/c,q\tilde t/b\big)M_{Nm}(\tilde t,k)D_{m }\big(\tilde t;b,c\big)\tilde{\alpha}_{m}.\label{RHS2}
\end{gather}
Equating expressions \eqref{LHS} and \eqref{RHS2}, cancelling common prefactors independent of the summation indices, and denoting $\alpha_{m}:=q^{-m(m+1)}\tilde{\alpha}_{m}$ we obtain the identity
\begin{gather*} %\label{MatrixBailey2}
M_{Nn}(a,k)D_{n}(a;b,c)M_{nm}\big(\tilde t,a\big)\alpha_{m}=D_{N}\big(k;q\tilde t/c,q\tilde t/b\big)M_{Nm}\big(\tilde t,k\big)D_{m}\big(\tilde t;b,c\big)\alpha_{m},
\end{gather*}
where we assume that there is a summation over repeated matrix indices. After removing arbitrary numbers $\alpha_m$ we see that this is precisely the key matrix Bailey lemma identity of~\cite{VS-IMRN-2002} as described above~\eqref{MatrixBailey}. The matrix relation~\eqref{MatrixBailey} was explicitly written in such a form in~\cite{VS-2008}. Our result consists in the demonstration that it is nothing else than the direct reduction of the integral star-triangle relation to the discrete form.

\section{Matrix realization of the Coxeter relations}

Let us rewrite now the matrix Bailey lemma identity \eqref{MatrixBailey} as a Coxeter relation of the permutation group $\mathbb{S}_3$. We take the three parameters set $\underline{t}=\big(\tilde{t},a,k\big)$ and introduce generators of elementary permutations~$s_j$:
\begin{gather*}
s_1\big(\tilde{t},a,k\big)=\big(a,\tilde{t},k\big), \qquad s_2\big(\tilde{t},a,k\big)=\big(\tilde{t},k,a\big).
\end{gather*}
They satisfy simple quadratic Coxeter relations $s_i^2=1$ and the cubic one $s_1s_2s_1=s_2s_1s_2$.

Now we define two operators $S_{1}$ and $S_{2}$ as the following matrices
\begin{gather*}
S_1\big(\tilde{t},a,k\big):=M\big(\tilde{t},a\big), \qquad S_2\big(\tilde{t},a,k\big):=D\big(\tilde{t};b,c\big),
\end{gather*}
where $b$ and $c$ are fixed in \eqref{bc}. These operators are assumed to use the twisted multiplication rule
\begin{gather*}
S_iS_j:=S_i(s_j\underline{t})S_j(\underline{t}),\qquad i,j=1,2.
\end{gather*}
Then one can easily check that
\begin{gather*}
 S_1^{2}=M\big(a,\tilde{t}\big)M\big(\tilde{t},a\big)=1,
\end{gather*}
which follows from the result of \cite{VS-IMRN-2002}. Permutation of $a$ and $k$ leads to the changes $b\to q\tilde t/c$ and $c\to q\tilde t/b$. As a result, we have
\begin{gather*}
 S_2^{2}=D\big(\tilde{t};q\tilde{t}/c,q\tilde{t}/b\big)D\big(\tilde{t};b,c\big)=1.
\end{gather*}

Consider now the following cubic combination of $S_j$-operators
\begin{gather*}
S_1S_2S_1:=S_1(s_2s_1\underline{t})S_2(s_1\underline{t})S_1(\underline{t}) =S_1\big(a,k,\tilde{t}\big)S_2\big(a,\tilde{t},k\big)S_1\big(\tilde{t},a,k\big).
\end{gather*}
The far right $S_1$-factor coincides with $M(\tilde{t},a)$. As the permutation of $a$ with $\tilde{t}$ does not change~$b$ and~$c$, we have $S_2\big(a,\tilde{t},k\big)=D(a;b,c)$. By the definition we have also $S_1\big(a,k,\tilde{t}\big)=M(a,k)$. So, $S_1S_2S_1$ coincides with the left-hand side of the key Bailey lemma equality~\eqref{MatrixBailey}:
\begin{gather*}
S_1S_2S_1=M(a,k)D(a;b,c) M\big(\tilde{t},a\big).
\end{gather*}

Now we consider another cubic combination of $S_j$-operators
\begin{gather*}
S_2S_1S_2:=S_2(s_1s_2\underline{t}) S_1(s_2\underline{t})S_2(\underline{t})=S_2\big(k,\tilde{t},a\big)S_1\big(\tilde{t},k,a\big)S_2\big(\tilde{t},a,k\big).
\end{gather*}
By definition, $S_2\big(\tilde{t},a,k\big)=D\big(\tilde{t};b,c\big)$ and $S_1\big(\tilde{t},k,a\big)=M\big(\tilde{t},k\big)$. The far left $S_2$-operator can be written in the form
\begin{gather*}
S_2\big(k,\tilde{t},a\big)=D\left(k;y\sqrt{\frac{pqk\tilde{t}}{a}},y^{-1}\sqrt{\frac{qk\tilde{t}}{pa}}\right)=D\left(k;\frac{q\tilde{t}}{c},\frac{q\tilde{t}}{b}\right).
\end{gather*}
Therefore we have
\begin{gather*}
S_2S_1S_2=D\big(k;q\tilde{t}/c,q\tilde{t}/b\big)M\big(\tilde{t},k\big)D\big(\tilde{t};b,c\big).
\end{gather*}
Finally, we see that the cubic Coxeter relation{\samepage
\begin{gather*}
S_1S_2S_1=S_2S_1S_2.
\end{gather*}
is nothing else than the key Bailey lemma identity \eqref{MatrixBailey}.}

To conclude, in this paper we have described direct reduction of the integral Bailey lemma to the matrix form. We interpreted also the key matrix identity of the latter lemma as the star-triangle relation, or the cubic Coxeter relation for permutation group $\mathbb{S}_3$.

As a next step of the development of derived results it is necessary to apply them to the problem of reducing the Yang--Baxter equation built with the help of integral Bailey lemma in~\cite{DS}. For that it is necessary to understand whether it is possible to extend the construction of $\mathbb{S}_3$-group generators given above to the group~$\mathbb{S}_4$ with appropriate interpretation of the corresponding Coxeter relations. An investigation of the relation of matrix Bailey lemma with the Sklyanin algebra is another problem for future considerations.

Another relevant subject concerns the supersymmetric field theories. Elliptic hypergeometric integrals are known to describe superconformal indices of such theories in four dimensions. The residue calculus corresponds in this picture to giving vacuum expectation values to some fields or insertion of surface defects and the results of the present paper may be useful for that interpretation as well. A comprehensive survey of this topic can be found in~\cite{RR}.

\subsection*{Acknowledgements}
This work is partially supported by Laboratory of Mirror Symmetry NRU HSE, RF government grant, ag.~no.~14.641.31.0001.

\pdfbookmark[1]{References}{ref}
\LastPageEnding

\end{document}